\documentclass[11pt, twoside, eqno]{article}
\usepackage{latexsym}
\usepackage{amssymb}
\usepackage{amsfonts}
\usepackage{graphicx}
\begin{document}
\begin{center}

\noindent {\bf \Large On completely regular and strongly regular ordered $\Gamma$-semigroups
}\bigskip

{\bf Niovi Kehayopulu}\bigskip

July 16, 2013
\end{center}
{\footnotetext{Department of Mathematics, University of Athens,
15784 Panepistimiopolis, Athens, Greece.

\hspace{0.2cm} email: nkehayop@math.uoa.gr.}
\bigskip

\noindent \bigskip

Our aim is to show the way we pass from the results of ordered semigroups (or semigroups) to ordered $\Gamma$-semigroups (or $\Gamma$-semigroups). The results of this note have been transferred from ordered semigroups. The concept of strongly regular $po$-$\Gamma$-semigroups has been first introduced here and a characterization of strongly regular $po$-$\Gamma$-semigroups is given. \bigskip

\noindent{\bf Definition 1}. Let $M$ be a $po$-$\Gamma$-semigroup. A
nonempty subset $B$ of $M$ is called a
{\it bi-ideal} of $M$ if the following assertions are satisfied:

(1) $B\Gamma M\Gamma B\subseteq B$ and

(2) if $a\in M$ and $b\in M$ such that $b\le a$, then $b\in
B$.\smallskip

For a subset $A$ of $S$ we denote by $B(A)$ the bi-ideal of $S$ generated by $A$ and we have $B(a)=(a\cup A\Gamma
A\Gamma A]$. \medskip

\noindent{\bf Proposition 2.} (cf. [2]) {\it Let M be a
$po$-$\Gamma$-semigroup and $B(x)$, $B(y)$ the bi-ideals of $M$
generated by the elements $x$ and $y$ of $M$, respectively, then
we have $$B(x)\Gamma M\Gamma B(y)\subseteq (x\Gamma M\Gamma y].$$}
{\bf Proposition 3.} (cf. [3]) {\it A
$po$-$\Gamma$-semigroup M is completely regular if and only if for
every $a\in M$ there exist $x\in M$ and $\gamma,\mu,\rho,
\xi\in\Gamma$ such that $$a\le (a\gamma a)\mu x\rho (a\xi
a).$$}{\bf Proof.} $\Longrightarrow$. By hypothesis, we
have$$a\le a\zeta t\omega a, \; a\le a\gamma a\mu y \mbox { and }
a\le z\rho a\xi a$$ for some $t,y,z\in M$ and
$\zeta,\omega,\gamma,\mu,\rho,\xi\in\Gamma$. Then we have

$$a\le a\zeta t\omega a\le (a\gamma a\mu y)\zeta t\omega (z\rho
a\xi a)=(a\gamma a)\mu (y\zeta t\omega z)\rho (a\xi a).$$We have
$(y\zeta t)\omega z\in M\Gamma M\subseteq M$. We put $x:=y\zeta
t\omega z$, and we have $a\le (a\gamma a)\mu x\rho (a\xi a)$,
where $x\in M$ and $\gamma,\mu,\rho,\xi\in\Gamma$.\smallskip

$\Longleftarrow$. Let $a\in M$. By hypothesis, there exist $x\in
M$ and $\gamma,\mu,\rho, \xi\in\Gamma$ such that $$a\le (a\gamma
a)\mu x\rho (a\xi a)\le a\gamma (a\mu x\rho a)\xi a,\; (a\gamma
a)\mu (x\rho a\xi a), \;(a\gamma a\mu x)\rho (a\xi a)$$where $a\mu
x\rho a$, $x\rho a\xi a$, $a\gamma a\mu x\in M$, so $M$ is
regular, left regular and right regular.\medskip

\noindent{\bf Proposition 4.} (cf. [3]) {\it A
$po$-$\Gamma$-semigroup M is completely regular if and only if
every bi-ideal B of M is semiprime}.$\hfill\Box$\medskip

\noindent {\bf Proof.} $\Longrightarrow$. Let $B$ be
a bi-ideal of $M$, $a\in M$ and $a\Gamma a\subseteq B$. Since $M$
is completely regular, by Proposition 3,
there exist $x\in M$ and $\gamma,\mu,\rho, \xi\in\Gamma$ such that
$$a\le (a\gamma a)\mu x\rho (a\xi a)\in (a\Gamma a)\Gamma M\Gamma
(a\Gamma a)\subseteq
B\Gamma M\Gamma B\subseteq B,$$so $a\in B$.\\
\noindent$\Longleftarrow$. Let $a\in M$. The nonempty set
$(a\Gamma a\Gamma M\Gamma a\Gamma a]$ is a bi-ideal of $M$. In
fact: Let $x,y\in (a\Gamma a\Gamma M\Gamma a\Gamma a]$, $\alpha,
\beta\in\Gamma$ and $z\in M$. We have $x\le a\gamma a\rho u\mu
a\xi a$ and $y\le a\zeta a\delta v\sigma a\lambda a$ for some
$u,v\in M$ and
$\gamma,\rho,\mu,\xi,\zeta,\delta,\sigma,\lambda\in\Gamma$. Then
we have\begin{eqnarray*}x\alpha z\beta y&\le& (a\gamma a\rho u\mu
a\xi a)\alpha z\beta (a\zeta a\delta v\sigma a\lambda
a)\\&=&a\gamma a\rho(u\mu a\xi a\alpha z\beta a\zeta a\delta
v)\sigma a\lambda a\in a\Gamma a\Gamma M\Gamma a\Gamma
a,\end{eqnarray*} so $x\alpha z\beta y\in (a\Gamma a\Gamma M\Gamma
a\Gamma a]$. Let now $x\in (a\Gamma a\Gamma M\Gamma a\Gamma a]$
and $y\in M$ such that $y\le x$. Then we have $y\in ((a\Gamma a\Gamma M\Gamma a\Gamma a]]=(a\Gamma a\Gamma M\Gamma a\Gamma a]$. The rest
of the proof is as in [1].\medskip

\noindent{\bf Proposition 5.} (cf. [4]) {Let M be a $po$-$\Gamma$-semigroup.
The following are equivalent:

$(1)$ M is completely regular.

$(2)$ $B(a)=B(a\Gamma a)=B(a\Gamma a\Gamma M\Gamma a\Gamma a)$
for every $a\in M$.

$(3)$ $B(a)=B(a\Gamma a)$.}\medskip

\noindent{\bf Proof.} $(1)\Longrightarrow (2)$. Let $a\in M$.
Since $M$ is regular, for the element $a$ of $M$, we have
$B(a)=(a\Gamma M\Gamma a]$ and, for the subset $a\Gamma a$ of $M$,
we have $B(a\Gamma a)=(a\Gamma a\Gamma M\Gamma a\Gamma a]$. Since
$M$ is right regular and left regular, we have
\begin{eqnarray*}(a\Gamma M\Gamma a]&\subseteq&((a\Gamma a\Gamma
M]\Gamma M\Gamma (M\Gamma a\Gamma a]]=((a\Gamma a\Gamma M]\Gamma
(M]\Gamma (M\Gamma a\Gamma a]]\\&\subseteq&((a\Gamma a\Gamma
M)\Gamma M\Gamma (M\Gamma a\Gamma a)]\subseteq (a\Gamma a\Gamma
M\Gamma a\Gamma a]\subseteq(a\Gamma M\Gamma a],\end{eqnarray*}so
$(a\Gamma M\Gamma a]=(a\Gamma a\Gamma M\Gamma a\Gamma a]$. Thus we
have $$B(a)=(a\Gamma M\Gamma a]=(a\Gamma a\Gamma M\Gamma a\Gamma
a]=B(a\Gamma a).$$In addition,\begin{eqnarray*}B(a\Gamma a\Gamma
M\Gamma a\Gamma a)&=& ((a\Gamma a\Gamma M\Gamma a\Gamma a)\cup
(a\Gamma a\Gamma M\Gamma
a\Gamma a)\Gamma M\Gamma (a\Gamma a\Gamma M\Gamma a\Gamma a)]\\
&=&(a\Gamma a\Gamma M\Gamma a\Gamma a]=B(a\Gamma
a).\end{eqnarray*}Thus we obtain $B(a)=B(a\Gamma a)=B(a\Gamma
a\Gamma M\Gamma a\Gamma a)$.\smallskip

\noindent $(3)\Longrightarrow (1)$. Let $a\in M$. By
hypothesis,$$a\in B(a)=B(a\Gamma a)=(a\Gamma a\cup a\Gamma a\Gamma
M\Gamma a\Gamma a]=(a\Gamma a]\cup (a\Gamma a\Gamma M\Gamma
a\Gamma a].$$If $a\le a\gamma a$ for some $\gamma\in \Gamma$, then
$a\le (a\gamma a)\gamma (a\gamma a)\le a\gamma a\gamma a\gamma
(a\gamma a)$ and, by Proposition 3, $M$ is completely
regular. If $a\in (a\Gamma a\Gamma M\Gamma a\Gamma a]$ then again
by Proposition 3, $M$ is completely regular.
\medskip

\noindent{\bf Proposition 6.} {\it Let $M$ be an $po$-$\Gamma$-semigroup. If $M$ is completely regular, then for each bi-ideal
$B$ of $M$, we have $B=(B\Gamma B]$. "Conversely", if $M$ has the
property $B=(B\Gamma B]$ for each bi-ideal $B$ of $M$, then $M$ is
regular. }\medskip

\noindent{\bf Proof.} $\Longrightarrow$. Let $B$ be a bi-ideal of
$M$. Then $B\Gamma M\Gamma B\subseteq B$. Since $M$ is regular, we
have $B\subseteq (B\Gamma M\Gamma B]$. Thus we have $B\subseteq
(B\Gamma M\Gamma B]\subseteq (B]=B$, and $B=(B\Gamma M\Gamma B]$.
Then we have\begin{eqnarray*}B\Gamma B&=&(B\Gamma M\Gamma B]\Gamma
(B]\subseteq (B\Gamma M\Gamma B\Gamma B]\\&\subseteq&(B\Gamma
M\Gamma B]=B.\end{eqnarray*} and $(B\Gamma B]\subseteq (B]=B$. On
the other hand, since $M$ is completely regular, we
have\begin{eqnarray*}B&\subseteq& (B\Gamma B\Gamma M\Gamma B\Gamma
B]\subseteq ((B\Gamma M\Gamma B)\Gamma B]\\&=&((B\Gamma M\Gamma
B]\Gamma B]=(B\Gamma B].\end{eqnarray*}Therefore we have
$B=(B\Gamma B]$.\smallskip

\noindent$\Longleftarrow$. Let $a\in M$. By hypothesis, we have
\begin{eqnarray*}a\in B(a)&=&(B(a)\Gamma B(a)]=((B(a)\Gamma B(a)]
\Gamma B(a)]\\&=&((B(a)\Gamma B(a))\Gamma B(a)\subseteq (B(a)\Gamma M\Gamma B(a)]\\
&=&((a\cup a\Gamma M\Gamma a]\Gamma M\Gamma (a\cup a\Gamma M\Gamma a]]\\
&=&((a\cup a\Gamma M\Gamma a]\Gamma (M]\Gamma (a\cup a\Gamma M\Gamma a]]\\
&=&((a\cup a\Gamma M\Gamma a)\Gamma M\Gamma (a\cup a\Gamma M\Gamma
a)]\\&=&(a\Gamma M\Gamma a\cup a\Gamma M\Gamma a\Gamma M\Gamma
a\cup a\Gamma M\Gamma a\Gamma M\Gamma a\Gamma M\Gamma
a]\\&=&(a\Gamma M\Gamma a],
\end{eqnarray*}so $S$ is regular.\medskip

\noindent{\bf Definition 7.} A $po$-$\Gamma$-semigroup $M$ is called
{\it strongly regular} if for every $a\in M$ there exist $x\in M$
and $\gamma, \mu\in \Gamma$ such that$$a\le a\gamma x\mu a \mbox {
and } a\gamma x=x\gamma a=x\mu a=a\mu x.$$We remark that if $M$ is
a strongly regular $po$-$\Gamma$-semigroup, then it is left
regular, right regular and regular. In fact: Let $a\in M$. Since
$M$ is strongly regular, there exist $x\in M$ and $\gamma,
\mu\in\Gamma$ such that $a\le a\gamma x\mu a$ and $a\gamma
x=x\gamma a=x\mu a=a\mu x$. Since $a\le (a\gamma x)\mu a=(x\gamma
a)\mu a=x\gamma a\mu a$, $M$ is left regular. Since $a\le a\gamma
(x\mu a)=a\gamma (a\mu x)=a\gamma a\mu x$, $M$ is right regular.
$M$ is clearly regular as well, so $M$ is completely regular.
\medskip

\noindent{\bf Theorem 8.} (cf. [5]) {\it Let M be a strongly regular
$po$-$\Gamma$-semigroup. Then, for every $a\in M$, there exist
$y\in M$ and $\gamma, \mu\in \Gamma$ such that$$a\le a\gamma y\mu
a,\; y\le y\mu a\gamma y \mbox { and } a\gamma y=y\gamma a=y\mu a=
a\mu y.$$}{\bf Proof.} Let $a\in M$. Since $M$ is strongly
regular, there exist $x\in M$ and $\gamma, \mu\in \Gamma$ such
that $a\le a\gamma x\mu a$ and $a\gamma x=x\gamma a=x\mu a=a\mu
x$. Then we have $$a\le a\gamma x\mu a\le (a\gamma x\mu a)\gamma
x\mu a=a\gamma (x\mu a\gamma x)\mu a.$$We put $y:=x\mu a\gamma x$,
and we have

$a\le a\gamma y\mu a$,\begin{eqnarray*}y=x\mu a\gamma
x&\le&x\mu(a\gamma x\mu a)\gamma x= (x\mu a\gamma x)\mu a\gamma
x=y\mu a\gamma x\le y\mu (a\gamma x\mu a)\gamma x\\&=&y\mu
a\gamma(x\mu a\gamma x)=y\mu a\gamma y,
\end{eqnarray*}so $y\le y\mu a\gamma y$,
\begin{eqnarray*} a\gamma y&=&a\gamma (x\mu a\gamma x)=(a\gamma
x)\mu (a\gamma x)=(x\mu a)\mu (x\gamma a)=x\mu (a\mu x)\gamma
a\\&=&x\mu (a\gamma x)\gamma a=(x\mu a\gamma x)\gamma a=y\gamma a,
\end{eqnarray*}$$y\gamma a=(x\mu a\gamma x)\gamma a=x\mu a\gamma (x\gamma a)=x\mu
a\gamma (x\mu a)=(x\mu a\gamma x)\mu a=y\mu a,$$
{\begin{eqnarray*}y\mu a&=&(x\mu a\gamma x)\mu a=x\mu (a\gamma
x)\mu a=x\mu (a\mu x)\mu a=(x\mu a)\mu (x\mu a)=(a\mu x)\mu
(a\gamma x)\\&=&a\mu (x\mu a\gamma x)=a\mu y.\end{eqnarray*}
\noindent{\bf Theorem 9.} (cf. [5]) {\it Let M be a $po$-$\Gamma$-semigroup. The following
are equivalent:}\begin{enumerate}
\item [(1)] {\it M is strongly regular}.
\item [(2)] {\it M is left regular, right regular, and $(M\Gamma a\Gamma M]$
is a strongly regular subsemigroup of M for every $a\in M$}.
\item [(3)] {\it For every $a\in M$, we have $a\in (M\Gamma a]\cap (a\Gamma
M]$, and $(M\Gamma a\Gamma M]$ is a strongly regular subsemigroup
of M.}\end{enumerate} }\noindent{\bf Proof.} $(1)\Longrightarrow (2)$. Let $a\in M$. The set $(M\Gamma
a\Gamma M]$ is a strongly regular subsemigroup of $M$. In fact:
Let $b\in (M\Gamma a\Gamma M]$. Since $b\in M$ and $M$ is strongly
regular, by Proposition 8, there exist $x\in M$ and $\gamma, \mu\in
\Gamma$ such that $b\le b\gamma x\mu b$, $x\le x\mu b\gamma x$,
and $b\gamma x=x\gamma b=x\mu b=b\mu x$. It is enough to prove
that $x\in (M\Gamma a\Gamma M]$.

Since $b\in (M\Gamma a\Gamma M]$, there exist $z,t\in M$ and
$\xi,\rho\in \Gamma$ such that $b\le z\xi a\rho t$. Then we have
$$x\le x\mu b\gamma x\le x\mu (z\xi a\rho t)\gamma x\in (M\Gamma
M)\Gamma a\Gamma (M\Gamma M)\subseteq M\Gamma a\Gamma M,$$ so
$x\in (M\Gamma a\Gamma M]$. Moreover, $M$ is left regular and
right regular, and (2) holds.\smallskip

\noindent$(2)\Longrightarrow (3)$. Let $a\in M$. Since $M$ is left
regular, we have $$a\in (M\Gamma a\Gamma a]\subseteq ((M\Gamma
M)\Gamma a]\subseteq (M\Gamma a].$$Since $M$ is right regular, we
have $a\in (a\Gamma a\Gamma M]\subseteq (a\Gamma (M\Gamma
M)]\subseteq (a\Gamma M]$, thus we get $a\in (M\Gamma a]\cap
(a\Gamma M]$, and condition (3) is satisfied.\smallskip

\noindent$(3)\Longrightarrow (1)$. Let $a\in M$. Since $a\in
(M\Gamma a]\cap (a\Gamma M]$, there exist $z,t\in M$ and $\xi,
\rho\in \Gamma$ such that $a\le z\xi a$ and $a\le a\rho t$. Then
we have$$a\le a\rho t\le (z\xi a)\rho t=z\xi (a\rho t)\in
M\Gamma a\Gamma M,$$ then $a\in (M\Gamma a\Gamma M]$. Since
$(M\Gamma a\Gamma M]$ is strongly regular, there exists $x\in
(M\Gamma a\Gamma M] (\subseteq M)$ and $\gamma, \mu\in\Gamma$ such
that $a\le a\gamma x\mu a$ and $a\gamma x=x\gamma a=x\mu a=a\mu x$, thus
$M$ is strongly regular.\bigskip

\noindent{\bf Conclusion.} If we want to get a result on a
$po$-$\Gamma$-semigroup, we first solve it for a $po$-semigroup
and then we have to be careful to define the analogous concepts in
case of the $po$-$\Gamma$-semigroup (if they do not defined
directly) and put the $"\Gamma"$ in the appropriate place. We
never solve the problem directly in $po$-$\Gamma$-semigroups.
{\small

\end{document}